\documentclass[10pt,a4paper]{article}

\usepackage{color}
\usepackage{ifpdf}
\ifpdf
 \usepackage[pdftex]{graphicx}
\else
 \usepackage[dvips]{graphicx}
\fi

\usepackage{hyperref}

\usepackage[T1]{fontenc}
\usepackage{a4,amssymb,amsthm,amsmath}
\usepackage{lmodern}
\usepackage[all]{xy}
\def\à{\`a}
\def\è{\`e}
\def\ä{\"a}

\newcommand{\g}{\mathfrak{h}}
\newcommand{\R}{\mathbb{R}}
\newcommand{\C}{\mathbb{C}}

\newcommand{\X}{\mathcal{X}}
\newcommand{\Y}{\mathcal{Y}}

\newcommand{\T}{\mathbb{T}}

\newcommand{\V}{\mathcal{V}}
\newcommand{\NN}{\mathcal{N}}
\newcommand{\HH}{\mathcal{H}}

\newcommand{\Sp}{\mathbb{S}}

\newcommand{\ZZ}{\mathcal{Z}}

\newcommand{\RR}{\mathcal{R}}
\newcommand{\TT}{\mathcal{T}}

\newcommand{\I}{\mathcal{I}}
\newcommand{\J}{\mathcal{J}}
\newcommand{\K}{\mathcal{K}}
\newcommand{\JJ}{\mathbb{J}}

\newcommand{\lra}{\longrightarrow}
\newcommand{\lms}{\longmapsto}

\newcommand{\w}{\wedge}
\newcommand{\wh}{\widehat}

\theoremstyle{definition}

\author{Guillaume~Deschamps}
 \title{   Twistor space of a generalized quaternionic manifold.}
 \date{\today}

\begin{document}

\maketitle

\begin{center}
{\bf Abstract}
\end{center}

We first make a little survey of the twistor theory for hypercomplex, generalized hypercomplex, quaternionic or generalized quaternionic manifolds.
This last theory was iniated by Pantilie \cite{Pan}, who shows that any generalized almost quaternionic manifold equipped with an appropriate connection admit a twistor space with an almost generalized complex structure.

 The aim of this article is to give an integrability criterion for this generalized almost complex structure and to give some examples  especially in the case of generalized hyperkähler manifolds using the generalized Bismut connection, introduced by Gualtieri \cite{Gua3}.

\setcounter{tocdepth}{4}
\tableofcontents
\section{Introduction}
In this article we review some properties of KT and HKT-manifolds, which closely resemble those of Kähler and hyperkähler ones respectively.  In physics, KT and HKT-manifolds arise as target spaces of two-dimensional supersymmetric sigma models with Wess-Zumino term \cite{GHR}. Another application of these geometries is in the context of black-holes, where the moduli spaces of a class of black-hole supergravity solutions are HKT-manifolds \cite{GPS}. Homogeneous manifolds have been investigated, and they have found many applications in physics in the context of sigma models and supergravity theory \cite{HW}. In mathematics, these notions are closely related with generalized geometry introduced by Hitchin \cite{Hit} and clarified by Gualtieri \cite{Gua2}. In section 3 and 5 we will recall these relations.

In section 6, we will make a little survey of twistor theory for hypercomplex and almost quaternionic  manifolds. The idea of a twistor space is to encode the geometric properties of the target manifold $M$ in term of holomorphic structure of $Z$.

\quad\\
{\bf Theorem 1 \cite{ AHS78,Sal1,Sal2,Bes87}.} Let $n\geq0$ and let $(M,Q)$ be an almost quaternionic $4n$-manifold. If $\nabla$ is a connection on $TM$ compatible with $Q$ then its twistor space admit a natural almost complex structure $\JJ_\nabla$ which is integrable if and only if, with respect to all local almost complex structures $J$ leaving in  $Q$ and all  sections $X,Y$ of $T M$ :
\begin{enumerate}
\item[(1)] The torsion $T$ of $\nabla$  satisfies :
$$
T(JX,JY)-JT(JX,Y)-JT(X,JY)-T(X,Y)=0.$$

\item[(2)] The curvature $R$ of $\nabla$  satisfies :
$$
\Big(R(X\w Y-JX\w JY)+JR(JX\w Y+X\w JY)\Big).J=0.
$$
\end{enumerate}

\quad\\
This theorem allows to give many examples of complex twistor space,  in particular, we have :

\quad\\
{\bf Theorem 2  \cite{Sal1,Sal2,Bes87}.} Let $n>1$ and let  $(M,Q)$ be an almost quaternionic $4n$-manifold. If $\nabla$ is a torsion free connection on $TM$ compatible with $Q$  then $\JJ_\nabla$ is integrable.

 \quad\\
 The purpose of this article is to extend these theorems in the context of generalized geometry. Indeed,   Pantilie \cite{Pan} noticed that we can still defined a twistor space for any generalized almost quaternionic manifold $(M,\mathcal Q)$, and when $M$ admit a connection $\nabla$ on $TM\oplus T^\star M$ compatible with $\mathcal Q$ then its twistor space admit a natural generalized almost complex structure $\mathbb J_\nabla$. Theorem A, anwers the question of the integrability of $\JJ_\nabla$ (see section 6 for precise definitions) :
 
\quad\\
{\bf Theorem A.} Let $n\geq 0$ and  $(M,\mathcal Q,\nabla)$ be a generalized almost quaternionic $4n$-manifold with a connection $\nabla$ on $\T M$ compatible with $\mathcal Q$.  The generalized almost complex structure $\JJ_\nabla$ on $\mathcal Z(\mathcal Q)$ is integrable if and only if with respect to all local generalized almost complex structures $u$  leaving in  $\mathcal Q$ and all  sections $\X,\Y,\ZZ$ of $\T M$:

\begin{enumerate}
\item[(C1)] The generalized torsion $\TT$ of $\nabla$  satisfies:
$$
\TT(\X,\Y,\ZZ)-\TT(\X,u\Y,u\ZZ)-\TT(u\X,\Y,u\ZZ)-\TT(u\X,u\Y,\ZZ)=0.$$
\item[(C2)] The generalized curvature $\RR$  of $\nabla$   satisfies:
$$\Big(\RR\big(\X\wedge\Y-u\X\wedge  u\Y\big)
+u\RR\big(u\X\wedge\Y+\X\wedge u\Y\big)\Big).u=0.$$
\end{enumerate}

\quad\\
This theorem enables us to give many new examples of generalized complex twistor space. In particular,

\quad\\
{\bf Theorem B.} Let $n\geq 0$. If $(M,G,\I,\J,\K)$ is a twisted generalized hyperkähler $4n$-manifold and  if $D$ is the generalized Bismut connection introduced by Gualtieri \cite{Gua3}, then the generalized almost complex structure $\JJ_D$ on its twistor space is integrable. 
 
\quad\\
This result is motivated by the fact that generalized hyperkähler structures appear in some branches of theoretical physics, such as string theory or in the context of certain supersymmetric sigma models \cite{GHR,HP,HP2,P}.

\section{KT-manifold}
 Let $(M,I,g)$ be a complex hermitian manifold and let $E\lra M$ be a fiber bundle. We denote by $\Gamma(E)$ the set of all smooth sections.
 A connection $\nabla : \Gamma(TM)\lra \Gamma(T^\star M\otimes TM)$ is called {\bf Hermitian} if $\nabla I=\nabla g=0$.
 Let $T(X,Y)=\nabla_XY-\nabla_YX-[X,Y]$ be the torsion tensor of type $(1,2)$. We denote by the same letter the torsion tensor of type $(0,3)$ given by $T(X,Y,Z)=g(X,T(Y,Z))$.
 
\quad\\
{\bf Definition.} A  hermitian connection is called  {\bf a Bismut connection} if $T$ is skew-symmetric.  The 3-form $T$ is then called the torsion form of the Bismut connection.

\quad\\
{\bf Proposition 1 \cite{FT}.} Let $(M,g,I)$ be a complex hermitian manifold and $w\in \wedge^{1,1}(X)$ the associated hermitian form. There exist a unique Bismut connection $\nabla^B$ and the torsion form is equal to $I dw$ that is :
$$
T(X,Y,Z)=dw(IX,IY,IZ).
$$
 If we denote by $\nabla^g$ the Levi-Civita connection of $g$, we have $\nabla^B=\nabla^g+\frac{1}{2}g^{-1}H$ that is :
$$
g(\nabla^B_XY,Z)=g(\nabla^g_XY,Z)+\frac{1}{2}T(X,Y,Z),
$$
for all vector fields $X,Y,Z$.

\quad\\
 Clearly if $dw=0$ then the Bismut connection is torsion-free and thus coincides with the Levi-Civita connection : the manifold $(M,g,I)$ is therefore Kähler.

Connection with skew-symmetric torsion play an important role in string physics. In the physics literature, a complex hermitian manifold $(M,g,I)$ with a Bismut connection is called a {\bf KT-manifold} (Kähler with torsion manifold). If in addition the torsion 3-form is closed then $(M,g,I)$ is said to be a {\bf strong KT-manifold}. By proposition 1, a manifold is therefore strong KT if and only if $\partial\bar\partial w=0$. For a complex surfaces, this is equivalent to Gauduchon metric. The strong KT-manifolds have been recently studied by many authors and they have also applications in type II string theory and in 2-dimensional supersymetric $\sigma$-models \cite{GHR,S,IP}. They also have relations with twisted generalized Kähler geometry as we are now going to see.

\section{Generalized complex structure}
\subsection{Courant bracket}
 Let $X,Y\in\Gamma(TM)$ be two vector fields  and $\xi,\eta\in\Gamma(T^\star M)$ be two 1-form. On  $\T M:=TM\oplus T^\star M$ there is an inner product : 
$$
<X+\xi,Y+\eta>=\frac{1}{2}\Big(\xi(Y)+\eta(X)\Big),
$$ 
and a Courant bracket, which is a skew-symmetric bracket defined  by 
$$
[X+\xi,Y+\eta]=[X,Y]+\mathcal L_X\eta-\mathcal L_Y\xi-\frac{1}{2}d(i_X\eta-i_Y\xi),
$$
 where $[X,Y]$ is the Lie bracket. The Courant bracket on $\T M$ can be twisted by a real closed 3-form $h$ defining another bracket \cite{Gua2, SW}
$$
[X+\xi,Y+\eta]_h=[X+\xi,Y+\eta]+i_Yi_Xh.
$$
In fact this bracket defined a Courant algebroïd structure on $\T M$.

\quad\\
When $b$ is a 2-form on $M$, we will denote by $e^b=\left(\begin{array}{cc}1&0\\b&1\end{array}\right)$ the transformation sending $X+\xi$ on $X+\xi+i_Xb$. This transformation is orthogonal for the inner product and is an automorphism for the Courant bracket if and only if $b$ is closed. 
\subsection{Generalized metric}
Let $M$ be a 2n-manifolds, since the bundle $\T M\lra M$ has a natural inner product, it has structure group $O(2n,2n)$.

\quad\\
{\bf Definition.}  A generalized metric is a reduction of the structure group from $O(2n,2n)$ to its maximal compact subgroup $O(2n)\times O(2n)$.

\quad\\
A generalized metric is equivalent to the choice of a 2n-dimensional subbundle $C^+$ which is positive definite with respect to the inner product. Let $C^-$ be the (negative definite) othogonal complement to $C^+$. Note that the splitting 
$$
\T M=C^+\oplus C^-$$ defines a positive definite metric on $\T M$ via :
$$
G=<.,.>\vert_{C^+}-<.,.>\vert_{C^-}.
$$
We denote by the same letter the isomorphism $G :\T M\lra \T M$ with $\pm 1$ eigenspace $C^\pm$, which is symmetric $G^\star =G$ and square to the identity $G^2=Id$.

\quad\\
{\bf Proposition 2 \cite{Gua2}.} A generalized metric is equivalent to specifying a riemannian metric $g$  and a 2-form $b$ on $TM$ such that :
\begin{enumerate}
\item[i)]
$G=e^b
\left(\begin{array}{cc}0&g^{-1}\\g&0\end{array}\right)
e^{-b}
$
\item[ii)] $C^\pm=\{X+(b\pm g)X\in\T M/X\in TM\}.$
\end{enumerate}

\subsection{Generalized complex structure}
A generalized almost complex structure on $M$ is an endomorphism $\J$ of $\T M$ which satisfies $\J^2=-1$ and $\J^\star=-\J$. That is a reduction of the structure group from $O(2n,2n)$ to $U(n,n)$.

\quad\\
{\bf Definition.} A generalized almost complex structure $\J$ is said to be a twisted generalized complex structure with respect to a closed 3-form $h$ when its $i$-eigenbundle $\T^{1,0}\subset \T M\otimes\C$ is involutive with respect to the $h$-twisted Courant bracket. We also said that $\J$ is twisted integrable or simply integrable when $h=0$.

\quad\\
Let $\NN_h$ be the Nijenhuis tensor of $\J$ defined on sections of $\T M$ by :
$$
\NN_h(\X,\Y)=[\J \X,\J\Y]_h-\J[\J\X,\Y]_h-\J[\X,\J\Y]_h-[\X,\Y]_h.
$$
When $h=0$ we simply note $\mathcal N$.

\quad\\
{\bf Proposition 3 \cite{Gua2}.} The twisted integrability of $\J$ is equivalent to the vanishing of the Nijenhuis tensor $\mathcal N_h$.

 \subsection{Generalized Kähler manifold}
Suppose that we have a generalized almost complex structure $\J$. To now reduce the structure group from $U(n,n)$ to $U(n)\times U(n)$ we need to choose a generalized metric $G$ which commutes with $\J$. Note that since $G^2=1$ and $G\J=\J G$, the map $G\J$ squares to $-1$ and since $G$ is symmetric and $\J$ is skew, $G\J$ is also skew, and therefore defines another generalized almost complex structure.

\quad\\
{\bf Definition \cite{Gua2}.} A reduction to $U(n)\times U(n)$ is equivalent to the existence of two commuting  generalized almost complex structures $\J_1$ and $\J_2$ such that $G=-\J_1\J_2$ is a generalized metric. We said that $(G,\J_1)$ is an almost generalized Kähler structure.

\quad\\
Since the bundle $C^+$ is positive definite while $TM$ is null, the projection $\pi : TM\oplus T^\star M\lra TM$ induces an isomorphism :
$$
\pi_\pm : C^\pm\lra TM.
$$
We denote by $P_{\pm}$ the projection from $\T M$ to $C^\pm$. Since $\J_1$ and $G$ commute, $\J_1$ stabilise $C^\pm$. By projection from $C^\pm$, $\J_1$ induces two almost complex structures on $TM$, which we denote $J^\pm$. They are compatible with the induced riemannian metric $g$ and the associated 2-forms are noted $w_\pm$. 
 Note that $\J_2=G\J_1$ implies that $\J_1=\J_2$ on $C^+$ and $\J_1=-\J_2$ on $C^-$.

 \quad\\
 {\bf Proposition 4 \cite{Gua2}.}  An almost generalized Kähler structure $(G,\J_1)$ is equivalent to the specification $(g,b,J^+,J^-)$ that is a riemannian metric $g$, a 2-form $b$ and two hermitian almost complex structures $J^\pm$ such that :

 \begin{enumerate}
\item[i)]
$G=e^b
\left(\begin{array}{cc}0&g^{-1}\\g&0\end{array}\right)
e^{-b}$
\item[ii)] $\J_1=\pi_+^{-1} J^+\pi P^++\pi_-^{-1} J^-\pi P^-$,

\item[iii)] $\J_2=\pi_+^{-1} J^+\pi P^+-\pi_-^{-1} J^-\pi P^-$,

\item[iv)] $\J_{1/2}=\displaystyle\frac{1}{2}e^b
\left(\begin{array}{cc}J_+\pm J_-&-(w_+^{-1}\mp w_-^{-1})\\w_+\mp w_-&-(J_+^\star\pm J_-^\star)\end{array}\right)
e^{-b}
$

\end{enumerate}

\quad\\
{\bf Definition \cite{Gua2}.} Let $(G,\J)$ be an almost generalized Kähler structure on $M$. When $\J$ and $G\J$ are both (twisted) generalized complex, we said that $(G,\J)$ is a (twisted) generalized Kähler structure on $M$.

\subsection{Relation betwenn KT and generalized Kähler manifold}
Let $(M,G,\J)$ be an almost generalized Kähler structure corresponding to the quadruple $(g,b,J_+,J_-)$.

\quad\\
{\bf Proposition 5 \cite{Gua2}.}  $(M,G, \J)$ is a twisted generalized Kähler structure if and only if:
\begin{enumerate}
\item[i)] $J_\pm$ integrable, and

\item[ii)] $h+db=-J_-dw_-=J_+dw_+$.

\end{enumerate}
This proposition shows that a twisted generalized Kähler structure on a riemannian manifold $(M,g)$  is the same that a bihermitian structure $(J_+,J_-)$ such that  the corresponding Bismut connections has torsions 3-forms which satisfy $T_+=-T_-$ and $dT_\pm=0$. In other words, a twisted generalized Kähler structure is a pair of strong KT-structures $(J_+,J_-)$ whose torsion satisfies $T_+=-T_-$

\quad\\
{\bf Proposition 6 \cite{Gua2}.} The torsion $T=-J_-dw_-=J_+dw_+$ of a twisted generalized Kähler structure is of type $(2,1)+(1,2)$ with respect to both complex structures $J_\pm$. Equivalently it satisfies the condition :
$$
T(X,Y,Z)-T(X,J_\pm Y,J_\pm Z)-T(J_\pm X,Y,J_\pm Z)-T(J_\pm X,J_\pm Y,Z)=0
$$
for all vector fields $X,Y,Z$ on $M$.

\quad\\
{\bf Example.} By \cite{SSTVP} any real compact Lie group $G$ of even dimension has a natural strong KT-metric and a twisted generalized Kähler structure \cite{Gua2}. 

\quad\\
Another notion due to physicists, is the notion of HKT-manifolds, which was suggested by Howe and Papadopoulos \cite{HP} and has been much studied since then.

\section{HKT-manifold}
{\bf Definition.} A riemannian $4n$-manifold $(M,g)$ admit a {\bf hypercomplex} structure if there exists a triple $(I,J,K)$ of complex structures such that  $IJ=-JI=K$. When each complex structure is compatible with the metric, we speack about {\bf hyperhermitian} structures.

\quad\\
Let $H\lra M$ be the vector bundle defined by $H=Vect(I,J,K)$. We said that a connection $\nabla$ on $TM$ is compatible with the hypercomplex structure or preserves the hypercomplex structure if $\nabla_X\sigma\in\Gamma(H)$ for all vector field $X$ and all smooth section $\sigma\in\Gamma(H)$.  On a hyperhermitian manifold, there are two natural torsion free connections, namely the Levi-Civita and the Obata connection. However, in general the Levi-Civita connection does not preserve the hypercomplex structure and the Obata connection does not preserve the metric. That's why we are interested in the following type of connections.

\quad\\
{\bf Definition.} A {\bf HKT-manifold} is a hyperhermitian $4n$-manifold $(M,g,I,J,K)$ with a connection $\nabla$ such that :
\begin{enumerate}
\item[i)] $\nabla g=0$

\item[ii)] $\nabla I=\nabla J=\nabla K=0$

\item[iii)] the torsion is totally skew-symmetric.
\end{enumerate}
 When the torsion is closed then $(M,g,I,J,K)$ is said to be a {\bf strong HKT-manifold}.

\quad\\
In contrast to the case of hermitian structure, not every hyperhermitian structure on a manifold admit a compatible HKT-connection but obviously, if such a connection exists, it is unique. Note that a HKT-connection is also the Bismut connection for each complex structure in the given hypercomplex structure. More generally we have

\quad\\
{\bf Proposition 7 \cite{CS}.} Let $(M,g,I,J,K)$ be an almost hyperhermitian manifold. It is an HKT-manifold if and only if :
\begin{eqnarray}
Idw_I=Jdw_J=Kdw_K
\end{eqnarray}
where $w_I, w_J,w_K$ are the associated hermitian form of $I,J$ and $K$. A holomorphic characterization has been given in \cite{GPP} where the autors proved that $(1)$ is equivalent to:
$$
\partial_I(w_J+iw_K)=0.
$$
 Many examples of HKT-manifolds have been obtained \cite{OP,GPP,V}. For instance, it has been shown that the geometry of the moduli space of a class of black holes in five dimensions is a HKT-manifold \cite{GPS}. 

\quad\\
We now consider (4,4)-supersymmetry structures on a riemannian manifold and see the link with HKT-structures. These structures were also introduced by Gates, Hull and Ro\u{c}ek \cite{GHR}, and formulated in Hitchin and Gualtieri's language as twisted generalized hyperkähler structures.

\section{Generalized hyperkähler manifold}
\subsection{Definition}
Let $(M,G)$ be a $4n$-manifold with a generalized metric. A (twisted) {\bf  generalized hyperhermitian} structure is a triple $(\I,\J,\K)$ of (twisted) generalized complex structures  such that 
\begin{enumerate}
\item[i)] $\I\J=-\J\I=\K$
\item[ii)]  $\I,\J,\K$ commute with $G$
\end{enumerate}

\quad\\
{\bf Definition.} A (twisted) {\bf generalized hyperkähler} structure on $M$ is a triple $(\I,\J,\K)$ of (twisted) generalized complex structure each of which forms a generalized Kähler structure with the same generalized metric $G$ and such that :
$$
\I\J=-\J\I=\K.
$$
{\bf Example 1 \cite{MV}.} A quaternionic Hopf surface $(\mathbb H-\{0\})/<q>$ where $q\in\R$, endowed with its two hypercomplex structures (left or right multiplication by $i,j$ or $k$) is an example of  a twisted generalized hyperkähler manifold.

\quad\\
{\bf Example 2 \cite{MV}.} Let $(M,g,I_\pm,J_\pm,K_\pm)$ be a generalized hyperkähler manifold of real dimension 4, and $E\lra M$ be a smooth complex vector bundle. Denote by $\mathcal M$ the moduli space of gauge-equivalence classes of anti-selfdual connections on $E$. Then $\mathcal M$ is equipped with a natural generalized hyperkähler structure.

\quad\\
{\bf Example 3 \cite{EG}.} The Neveu-Schwarz 5-branes solution provides an explicit example of generalized hyperkähler manifold found in string theory.

\subsection{Relation betwenn HKT and generalized hyperkähler structure}
Let $(M,G,\I,\J, \K))$ be an almost generalized hyperkähler structure corresponding to $(g,b,I_+,J_+,K_,I_-,J_-,K_-)$.

\quad\\
{\bf Proposition 8.}  $(M,G,\I, \J,\K)$ is a twisted generalized hyperkähler structure if and only if:
\begin{enumerate}
\item[i)] $(I_\pm, J_\pm,K_\pm)$ is a pair of hyperhermitian complex structure on $(M,g)$, and

\item[ii)] $\begin{array}{rrrrr}h+db&=&I_+dw_{I_+}=&J_+dw_{J_+}=&K_+dw_{K_+}\\
&=&-I_-dw_{I_-}=&-J_-dw_{J_-}=&-K_-w_{K_-}.\end{array}$

\end{enumerate}
In other words, a twisted generalized hyperkähler structure is a pair of strong HKT-structure $(I_\pm,J_\pm,K_\pm)$ whose torsion satisfies $T_+=-T_-$.

\quad\\
{\bf Corollary.} Torsions $T_+=-T_-$ of a twisted generalized hyperkähler structure is of type $(2,1)+(1,2)$ with respect to each complex structure $I_\pm,J_\pm$ or $K_\pm$.

\section{Twistor space}
In this section, we define the twistor space of a twisted generalized hyperkähler manifold $(M,G,\I,\J,\K)$ and more generally to a generalized almost quaternionic manifold. Unlike the approach of Bredhauer \cite{Bre}, still generalized on \cite{Des1}; in this paper, the twistor space is not an $\Sp^2\times\Sp^2$-fiber bundle but a $\Sp^2$-bundle exactly has in the original idea of Penrose \cite{Pen} and Salamon \cite{Sal1,Sal2}. We first review the results for a quaternionic manifold. 

\subsection{Twistor space of a quaternionic manifold}
Let $\Big(M,(I,J,K)\Big)$ be an hypercomplex $4n$-manifold, a triple of such  complex structures induces  a 2-sphere of integrable complex structures :
$$\{aI+bJ+cK/a^2+b^2+c^2=1\}.
$$
So it is natural to define the twistor space associated to this hypercomplex structure by:
$$Z=M\times\Sp^2=\{(m,aI+bJ+cK)/(a,b,c) \in\Sp^2\}.
$$
 The idea of a twistor space is to encode the geometric properties of the target manifold $M$ in the holomorphic structure of $Z$. Indeed, we are now going to define a natural almost complex structure $\JJ_\nabla$ or simply $\JJ$ for any connection $\nabla$ on $M$ that preserves the hypercomplex structure.
Such a connection induces decomposition:
$$
TZ=H\oplus  V
$$
of the tangent space of $Z$ into its vertical and horizontal component.
Since $\Sp^2$ has the natural complex structure of $\C P^1$, we take $\JJ\vert_V =\JJ_{\C P^1}$. But $H$ and $T M$ are isomorphic, so we may define $\JJ\vert_H$ by letting $\JJ\vert_H$ act at $(m,u)\in Z$ like $u$ on $T_mM$.\\

 The construction of a twistor space and its almost complex structure can be easily extended to any {\bf almost quaternionic} $4n$-manifold $(M,Q)$, that is manifold with a rank three subbundle $Q\subset End(TM)\lra M$ which is locally spanned by an almost hypercomplex structure $(I,J,K)$.  Such a locally defined triple $(I,J,K)$ will be called an admissible basis of $Q$.  A consequence of the definition of an almost quaternionic manifold is that the bundle $Q$ has structure group $SO(3)$. We then have a natural inner product on $Q$ by taking each admissible basis $(I,J,K)$ to be an orthonormal basis. The twistor space $Z(Q)$ of $(M,Q)$ is defined to be the unit sphere bundle of $Q$. This is a locally trivial bundle over $M$ with fiber $\Sp^2$. A linear connection $\nabla$ on $TM$ preserves $Q$ means that $\nabla_X\sigma\in\Gamma(Q)$ for all vector field $X$  and smooth section $\sigma\in\Gamma(Q)$. In this case, the same construction as before gives us an almost complex structure $\JJ_\nabla$ on $Z(Q)$  whose integrability depends on the torsion $T$ and the curvature $R$ of $\nabla$, defined by :
$$
\begin{array}{l}
T(X,Y)=\nabla_XY-\nabla_YX-[X,Y], \\
R(X,Y)=\nabla_X\nabla_Y-\nabla_Y\nabla_X-\nabla_{[X,Y]}.
\end{array}
$$
By skew symmetry we note $R(X\w Y)$ rather than $R(X,Y)$.

\quad\\
{\bf Theorem 1 \cite{ AHS78,Sal1,Sal2,Bes87}.} Let $n\geq0$ and let $(M,Q)$ be an almost quaternionic $4n$-manifold. If $\nabla$ is a connection on $TM$ compatible with $Q$ then its twistor space admit a natural almost complex structure $\JJ_\nabla$ which is integrable if and only if, with respect to all local almost complex structures $J$ leaving in  $Q$ and all  sections $X,Y$ of $T M$ :
\begin{enumerate}
\item[(1)] The torsion $T$ of $\nabla$  satisfies :
$$
T(JX,JY)-JT(JX,Y)-JT(X,JY)-T(X,Y)=0.$$

\item[(2)] The curvature $R$ of $\nabla$  satisfies :
$$
\Big(R(X\w Y-JX\w JY)+JR(JX\w Y+X\w JY)\Big).J=0.
$$
\end{enumerate}

\quad\\ In the particular case of a torsion free connection we have,

\quad\\
{\bf Theorem 2 \cite{Sal1,Sal2,Bes87}.} Let $n>1$ and let $(M,Q)$ be an almost quaternionic $4n$-manifold. If $\nabla$ is a torsion free connection on $TM$ compatible with $Q$, then $\JJ_\nabla$ is a complex structure on $Z(Q)$.

 \quad\\
 Pantilie \cite{Pan} extended this construction in the context of generalized geometry as we are now going to see.

 \subsection{Twistor space of a generalized hypercomplex manifold}

 For a generalized hypercomplex manifold $\Big(M,(\I,\J,\K)\Big)$ we can still defined the associated twistor space by :
$$\ZZ=M\times\Sp^2=\{(m,a\I+b\J+c\K)/(a,b,c)\in\Sp^2\},
$$
and we denote by $\pi_\ZZ : \ZZ\lra M$, the first projection.

\quad\\
As in the classical case, any connection on $\T M$ compatible with  the generalized hypercomplex structure defined  a natural generalized almost complex structure $\JJ_\nabla$ on $\ZZ$ and the construction is substantially the same. Indeed, since the connection allows us to parallel transport element $u$ of $\ZZ$, there is an associated decomposition $T\ZZ=\HH\oplus \mathcal V$ of the tangent space of $\ZZ$ into horizontal and vertical  parts. If we note $\HH^\star$ (resp. $\mathcal V^\star$) the element of $T^\star \ZZ$ null on $\mathcal V$ (resp. on $\HH$), it gives us a decomposition $\T \ZZ=(\HH\oplus\HH^\star)\oplus(\mathcal V\oplus\mathcal V^\star)$. The function $d\pi_\ZZ\oplus d\pi_\ZZ^\star$ gives us an isomorphism betwenn $\HH\oplus\HH^\star$ and $\T M$ which preserves the inner product. So we may define a bundle endomorphism :
$$
\JJ\vert_{\HH\oplus\HH^\star} :\HH\oplus\HH^\star\lra\HH\oplus\HH^\star,\quad \JJ\vert_{\HH\oplus\HH^\star}^2=-1
 $$
 by setting $\JJ\vert_{\HH\oplus\HH^\star}$ act at $(m,u)\in\ZZ$ by    $\Big(d\pi_\ZZ\oplus d\pi_\ZZ^\star\Big)^\star u$.
 On the other hand, $\mathcal V$ is just the tangent space to the fibers and so admits the natural  complex structure $\JJ\vert_{\mathcal V\oplus\mathcal{V}^\star}$ of $\C P^1$. This gives us a natural generalized almost complex structure $\JJ_\nabla$ on $\ZZ$ namely $\JJ\vert_{\HH\oplus\HH^\star}\oplus\JJ\vert_{\mathcal V\oplus\mathcal V^\star}$.

\subsection{Twistor space of a generalized quaternionic manifold}
We say that $M$ admit an {\bf generalized almost quaternionic} structure if there exists $\mathcal Q\lra M$  a rank three vector bundle $\mathcal Q\subset End(\T M)$ which is locally spanned by an almost generalized  hypercomplex structure.  The twistor space $\ZZ(\mathcal Q)$ of $(M,\mathcal Q)$ is still defined to be the unit sphere bundle of $\mathcal Q$ for the natural inner product in $\mathcal Q$ such that any admissible basis is othonormal. The bundle $\pi_{\ZZ(\mathcal Q)} : \ZZ(\mathcal Q)\lra M$  is a locally trivial bundle  with fibre $\Sp^2$ and structure group $SO(3)$.
Moreover it is not difficult to see that the former construction of the generalized almost complex structure $\JJ_\nabla$ associated to any connection $\nabla$ on $\T M$ preserving  $\mathcal Q$ works yet.

\quad\\
{\bf Extension.}  To simplify the notations, we extend $\nabla$ to $\T M$ asking $\nabla_\X=\nabla_{\pi(\X)}$ for all $\X\in\Gamma(\T M)$.

\quad\\ 
In his article Pantilie does not study the integrability of $\JJ$. In the next section we will give a criterion of integrability for $\JJ$. As in the usual case, it depend on  the generalized torsion $\TT$ and the generalized curvature $\RR$  of the connection $\nabla$. Recall that the  generalized torsion $\TT$, is defined by Gualtieri \cite{Gua3}, for all $\X,\Y,\ZZ\in\Gamma(\T M)$, by
 $$
 \TT(\X,\Y,\ZZ)=<\nabla_\X\Y-\nabla_\Y\X-[\X,\Y],\ZZ>+\frac{1}{2}\Big(<\nabla_\ZZ\X,\Y>-<\nabla_\ZZ\Y,\X>\Big).$$
 As $\nabla$ preserves the inner product, then $\TT$ is totally skew.
On the other hand, the generalized curvature  is defined by
 $$
 \RR(\X,\Y)=\nabla_\X\nabla_\Y-\nabla_\Y\nabla_\X-\nabla_{[\X,\Y]}.
 $$
  By skew symmetry, sometimes we will use the notation $\RR(\X\wedge\Y)$ rather than $\RR(\X,\Y)$. Since $\nabla$ is  an usual connection on $\T M$ then  $\RR$ is tensorial.
    
 \subsection{Integrability of the generalized almost complex structure}
 Let $(M,\mathcal Q)$ be a generalized almost quaternionic manifold. The data of a generalized almost hypercomplex structure $(\I,\J,\K)$ on an open set $\mathcal U$ of $M$ defines a trivialisation $\pi^{-1}_{\ZZ(\mathcal Q)}(\mathcal U)\simeq\mathcal U\times\Sp^2$. The local coordinates of a point in $\ZZ(\mathcal Q)$ will be denoted by $(m,u)$. 
The central theorem of this article is the following.

\quad\\
{\bf Theorem A.} Let $n\geq 0$ and  $(M,\mathcal Q,\nabla)$ be a generalized almost quaternionic $4n$-manifold with a connection $\nabla$ on $\T M$ compatible with $\mathcal Q$.  The  generalized almost complex structure $\JJ_\nabla$ on the twistor space $\ZZ(\mathcal Q)$ is integrable if and only if with respect to all local generalized almost complex structures leaving in  $\mathcal Q$, the two following conditions are satisfied:

\begin{enumerate}
\item[(C1)] The torsion $\TT$ is of type $(2,1)+(1,2)$. Equivalently it satisfies the local condition:
$$
\TT(\X,\Y,\ZZ)-\TT(\X,u\Y,u\ZZ)-\TT(u\X,\Y,u\ZZ)-\TT(u\X,u\Y,\ZZ)=0$$
for all sections $\X,\Y,\ZZ$ of $\T M$ and all $u=a \I+b\J+c\K$ with $(a,b,c)\in\Sp^2$.
\item[(C2)] The curvature form $\RR(\X,\Y,\ZZ,\mathcal U):=<\RR(\X,\Y)\ZZ,\mathcal U>$ has no component of type $(4,0)+(0,4)$. Equivalently it satisfies the local condition:
$$\Big(\RR\big(\X\wedge\Y-u\X\wedge  u\Y\big)
+u\RR\big(u\X\wedge\Y+\X\wedge u\Y\big)\Big).u=0$$
for all sections $\X,\Y$ of $\T M$ and all $u=a \I+b\J+c\K$ with $(a,b,c)\in\Sp^2$.
\end{enumerate}

\subsection{Proof}
We will used the notation $\wh\X\in\HH\oplus\HH^\star$ to denote the horizontal lift of a local smooth section of $\T M$, and we will speack about basic sections of $\HH\oplus\HH^\star$.

 Because the connection $\nabla$ preserve $\mathcal Q$, for all smooth sections $\X,\Y\in \T M$, and for all $u\in\Gamma(\mathcal Q)$, we have that $\RR(\X,\Y).u\in\Gamma(\mathcal Q)$ and more precisely it is a vertical element.
 In order to prove the integrability of $\JJ$, we consider its Nijenhuis tensor $\NN$, and the various components of it in the splitting  $\T \ZZ=(\HH\oplus\HH^\star)\oplus(\mathcal V\oplus\mathcal V^\star)$.
The computation of these components require the following proposition.

\quad\\
{\bf Proposition 9.} For all vertical vector fieds $A,B\in\Gamma(\mathcal V)$ and all horizontal lift $\wh{X+\xi}\in\Gamma(\HH \oplus \HH^\star)$ one has:

\begin{enumerate}
\item[i)] $[A,B]\in\Gamma(\mathcal V)$

\item[ii)] $[\wh X,A]\in \Gamma(\V)$

\item[iii)] $[\wh X+\wh \xi ,\JJ A]=\JJ[\wh X+\wh \xi,A]$

\item[iv)] $[\JJ(\wh X+\wh\xi),\JJ A]=\JJ[\JJ(\wh X+\wh \xi),A]$

\end{enumerate}

\quad\\
{\bf Proof.} The first point is a general fact for any vertical distribution, similarly the  second point is always true for any  basic vector field $\wh X$ and  any vertical vector field $A$ \cite{Bes87}.

The third formula follows from the parallel transport along horizontal directions respect the canonical metric and the orientation of the fibres, hence the vertical complex structure, so $[\wh X,\JJ A]=\JJ[\wh X,A]$. It remains to check that $[\wh \xi,A]=0=[\wh \xi,\JJ A]$ which is an immediate consequence of  the definition of the Courant bracket.

For the last formula, pick a local basis $(\X_1,\ldots,\X_{8n})$ of $\Gamma(\T M)$ and note $[\JJ_{ij}]$ the matrix of $\JJ$ in the basis $(\wh\X_1,\ldots,\wh\X_{8n})$ of $\HH\oplus\HH^\star$. Properties of the bracket  and the third point give us :

$$
\begin{array}{cccc}
&\JJ[\JJ\wh\X_j,A]&=&\JJ[\JJ_{ij}\wh\X_i,A]\\
&&=&\JJ\big(\JJ_{ij}[\wh\X_i,A]-A\wh\X_j\big)\\
&&=&\JJ_{ij}[\wh\X_i,\JJ A]-\JJ A\wh\X_j\\
\\\textrm{et}&[\JJ\wh\X_j,\JJ A]&=&[\JJ_{ij}\wh\X_i,\JJ A]\\
&&=&\JJ_{ij}[\wh\X_i,\JJ A]-\JJ A\wh\X_j.\; \square
\end{array}
$$

\quad\\
{\bf Corollary 1.} The component $\NN(\X,A)$ of the Nijenhuis tensor of $\JJ$ is null for all $\X\in \Gamma(\HH\oplus\HH^\star)$ and all $A\in \Gamma(\V)$.

\quad\\
{\bf Proof.} By linearity, one can supposed that $\X$ is  basic, the corollary is then an immediate consequence of points 3 and 4 of proposition 9. $\square$

\quad\\
{\bf Proposition 10.}  Let $\X,\Y\in\Gamma(\T M)$ be two local sections.  According to horizontal and vertical directions at a point $p=(m,u)\in\ZZ(\mathcal Q)$, one has:
$$[\wh\X,\wh\Y]=\wh{[\X,\Y]}+\RR(\X,\Y).u.
$$

\quad\\
{\bf Proof.} Identify $\mathbb H^n=\R^{4n}$ and let $\R^{4n\star}$ be the dual of $\R^{4n}$. The group $GL(2n,\mathbb H)$ (resp. $Sp(1)$) act on the right (resp. on the left) on $\R^{4n}\oplus\R^{4n\star}$. We note $\mathcal Sp(2n)$ the subgroup :
$$\mathcal S p(2n)=GL(2n,\mathbb H)\cap O(2n,2n).
$$
Let $G$ be the product $\mathcal Sp(2n)\mathcal Sp(1)$ and  $\mathcal P\lra M$ be the $G$-principal bundle. The twistor space $\ZZ(\mathcal Q)$ can be considered has
the associated fiber bunble of $\mathcal P$ with standard fiber $\Sp^2$. More precisely, the group $G$ acts on the right on
$\mathcal P\times\Sp^2$ by:
$$
\begin{array}{ccc}
\mathcal P\times \Sp^2\times G&\lra&\mathcal P\times\Sp^2\\
(q,j,g)&\lms&(q.g,g^{-1}.j)=(q.g,gjg^{-1})
\end{array}
$$
and $\ZZ(\mathcal Q)$ is the quotient of $\mathcal P\times\Sp^2$ by $G$. Note $\Pi$ the projection
$$
\begin{array}{cccc}
\Pi &: \mathcal P\times\Sp^2&\lra&\ZZ(\mathcal Q)\\
&(q,j)&\lms&u=q^{-1}jq.
\end{array}
$$
 Start by looking at the case where $ X,Y$ are two vector fields on  $M$.  As $\wh X,\wh Y$  are basic, then  \cite{Bes87} the horizontal part of $[\wh X,\wh Y]$ is precisely $\wh{[X,Y]}$. We denote by $\theta$ be the $G$-connection on $\mathcal P$, $ker \theta$ the associated horizontal distribution and
  $\widetilde X,\widetilde Y$ the horizontal lift of $X, Y$ in $\mathcal P$.  The vertical part of $[\widetilde X,\widetilde Y]$
is given by \cite{O},
\cite{Bes87} :
$$
(\theta\vert_{\V})^{-1}(\RR(X,Y)),$$
where by definition, $(\theta\vert_{\V})^{-1}(\RR(X,Y))$ is the vertical field on  $\mathcal P$ defined at the point $q\in \mathcal P$ by:
 $$
\frac{d}{dt}\vert_{t=0}\Big(q.\,exp(t\RR(X,Y))\Big)=p.\,\RR(X,Y).
$$
At $q\in\mathcal P$ we have
$d\Pi\Big(q.\RR(X,Y)\Big)=u\RR(X,Y)-\RR(X,Y)u=\RR(X,Y).u$. Then at the point $(m,u)\in\ZZ(\mathcal Q)$ we have 
$$
[\wh X,\wh Y]=\wh{[X,Y]}+\RR(X,Y).u
$$
Now let $X\in TM$ be a vector field and $\eta\in T^\star M$  be a 1-form on $M$. Using the definition of the Courant bracket we see that $ [\wh X,\wh \xi]=\wh{[X,\xi]}$ and from $\nabla_\xi=0$ we deduce that $\RR(X,\xi)=0$ and so :
$$
\begin{array}{rcl}
[\wh X,\wh \xi]&=&\wh{[X,\xi]}+\RR(X,\xi).u \\
0=[\wh\xi,\wh\eta]&=&\wh{[\xi,\eta]}+\RR(\xi,\eta).u=0\quad \square
\end{array}
$$

\quad\\
{\bf Corollary 2.} For all vertical 1-form $U^\sharp$ and all $\X\in\Gamma(\HH\oplus\HH^\star)$, the component $\NN(U^\sharp,\X)$ of the Nijenhuis tensor of $\JJ$ is  the horizontal form defined for all $\Y\in \Gamma(\HH)$ by :
$$
 <\NN(U^\sharp,\X),\Y>=
 U^\sharp\Big(\big(\RR(\X\wedge\Y-u\X\wedge u\Y)+
     u\RR(u\X\wedge\Y+\X\wedge u\Y)\big).u\Big).
 $$

 \quad\\
 {\bf Proof.} We denote by $\overrightarrow \X,\overrightarrow \Y$ the projection of $\X,\Y\in\HH\oplus \HH^\star$ over $\HH$. Using the definition of the Courant bracket we know that 
 $[U^\sharp,\X]=[U^\sharp,\overrightarrow{\X}]$ is a 1-form. More precisely for two vector fields $A\in \Gamma(\V)$ and $\overrightarrow{\Y}\in  \Gamma(\HH)$, at the point $p=(m,u)\in\ZZ(\mathcal Q)$ we have : 
$$
\begin{array}{ccc}
[U^\sharp,\X](A+\overrightarrow{\Y})&=&dU^\sharp(\overrightarrow{\X},\overrightarrow{\Y}+A)\\
&=&\overrightarrow{\X}.U^\sharp(A)-U^\sharp([\overrightarrow{\X},\overrightarrow{\Y}+A])\\
&=&\overrightarrow{\X}.U^\sharp(A)-U^\sharp([\overrightarrow{\X},A])
-U^\sharp\Big(\RR(\X\w\Y).u\Big)
\end{array}
$$
The point 3 of proposition 9 gives us $[\JJ
U^\sharp,\X](A)=\JJ[U^\sharp,\X](A)$, and so $\NN(U^\sharp,\X)$ is the horizontal 1-form defined by
$$
<\NN(U^\sharp,\X),\Y>=U^\sharp\Big(\big(\RR(\X\wedge\Y-u\X\wedge u\Y)+
     u\RR(u\X\wedge\Y+\X\wedge u\Y)\big).u\Big)
  \;\square
 $$

\quad\\
{\bf Corollary 3.} For all basic sections  $\wh\X,\wh\Y\in \Gamma(\HH\oplus\HH^\star)$, at the point $p=(m,u)\in\ZZ(\mathcal Q)$ :

\begin{enumerate}
\item[i)] the vertical part of $\NN(\wh\X,\wh\Y)$ is the following vector field :
$$-\Big(\RR\big(\X\wedge\Y-u\X\wedge u\Y\big)+
     u\RR\big(u\X\wedge\Y+\X\wedge u\Y\big)\Big).u$$
     
\item[ii)] the horizontal part of $\NN(\wh\X,\wh\Y)$ is a 1-form defined for any $\wh\ZZ\in\Gamma(\HH)$ by :
$$
<\NN(\wh\X,\wh\Y),\wh\ZZ>=\TT(\X,\Y,\ZZ)-\TT(\X,u\Y,u\ZZ)-\TT(u\X,\Y,u\ZZ)-\TT(u\X,u\Y,\ZZ)$$
\end{enumerate}

\quad \\
{\bf Proof.}     Pick an orthonormal basis $(\X_1, \ldots, \X_ {8n}) $ of $\T M$ defined above $ \mathcal U $. The distribution $ \mathcal
H \oplus \HH ^\star $ is stable for $ \JJ$, we note $ [\JJ_{ij}] $ his
matrix in the basic $ (\wh {\X_1}, \ldots, \wh {\X_{8n}})$.
By definition
$$\begin{array}{rll}
  \left[\mathbb J \widehat{\X_i},\mathbb
  J\widehat{\X_j}\right]&=&\overrightarrow{\JJ\wh{\X_i}}.(\JJ_{rj})\;\widehat{\X_r}
  -\overrightarrow{\JJ\wh{\X_j}}.(\JJ_{li})
\;\widehat{\X_l}+  \JJ_{li}\JJ_{rj}\left[\widehat{\X_l},\widehat{\X_r}\right]\\
&&-\JJ_{ri}d\JJ_{rj}+\JJ_{lj}d\JJ_{li}\\
  \,\left[\mathbb J\widehat{\X_i},\widehat{\X_j}\right]
  +\left[\widehat{\X_i},\mathbb
  J\widehat{\X_j}\right]  &=&-\overrightarrow{\wh{\X_j}}.(
  \JJ_{li})\;\widehat{\X_l}+\JJ_{li}\left[\widehat{\X_l},
  \widehat{\X_j}\right]+ \overrightarrow{\wh{\X_i}}.(
  \JJ_{rj})\;(\widehat{\X_r})
  +\JJ_{rj}\left[\widehat{\X_i},\widehat{\X_r}\right]\\
  &&+d\JJ_{ji}-d\JJ_{ij}.
\end{array}
$$
Using proposition 10, we deduce that at the point $p=(m,u)$ , the
vertical part of $\NN (\wh\X_i, \wh \X_j)$ is
$$
-\Big(\RR(\X_i\wedge\X_j-u\X_i\wedge u\X_j)+
     u\RR(u\X_i\wedge\X_j+\X_i\wedge u\X_j)\Big). u
$$
For the horizontal part, we consider a local section $s$ of
$ \ZZ (\mathcal Q) \lra M $ on $U$, such that $s(m) = u$ and $ (\nabla s)_m = 0$. This gives a local generalized almost complex structure $S$ on $U$. The
horizontal part of $ \NN (\wh \X_i, \wh \X_j) $ restricted to  $ s(M) $ is equal to the horizontal lift of the Nijenhuis tensor of $S$. Since the connection $\nabla $ has torsion $\TT$ and since $\nabla s=0$ at $m$ we have, at the point $(m,u)$ :
$$\begin{array}{llc}<\NN(\wh\X_i,\wh\X_j),\wh\X_k>&=&\TT(\X_i,\X_j,\X_k)-\TT(u\X_i,u\X_j, \X_k)\\
&& -\TT(u\X_i,\X_j,u\X_k)-\TT(\X_i,u\X_j,u\X_k).
 \end{array}$$

\quad\\
{\bf Proof of theorem A.}  Since fibers of $\ZZ(\mathcal Q)$ has the complex structure of $\C P^1$, we get $\NN(\mathcal U,\mathcal V)=0$ for all $\mathcal U,\mathcal V$ basic sections of $\mathcal V\oplus\mathcal V^\star$. The proof of theorem A is then an immediate consequence of the corollaries 1, 2 and 3.

\section{Applications}
\subsection{Generalized Bismut connection}
Generalized hyperkähler structures are among the simplest examples of generalized quaternionic structure $\mathcal Q$. In that case there is a natural connection preserving $\mathcal Q$. this connection was
introduced by Gualtieri \cite{Gua3} and is called generalized Bismut connection. We start by recalling its construction.
 Let $G=(g,b)$ be a generalized metric   and $C^+$ the associated maximal-positive-definite subbundle of $\T M$. Let $C : \T M\lra \T M$ be the automorphism  defined by $C(X+\xi)=X-\xi$. Write $\X=\X^++\X^-$ for the orthogonal projection of $\X\in\Gamma(\T M)$ to $C^\pm$ and let $h$ be any closed 3-form on $M$. 

\quad\\
{\bf Proposition 11 \cite{Gua3}.}   The operator :
$$
D_\X\Y=[\X^-,\Y^+]_{h}^++[\X^+,\Y^-]_{h}^-+[C\X^-,\Y^-]_{h}^-+[C\X^+,\Y^+]_{h}^+
$$
defines a connection on $\T M$, preserving both the inner product $<.,.>$ and the positive-definite metric $G$. So $D$ preserved $C^\pm$ and if we note $D^\pm$ the restriction of $D$ to $C^\pm$, then we have :
$$
D^\pm=\pi^{-1}_\pm\nabla^\pm\pi_\pm
$$
where $\nabla^\pm$ are the Bismut connection on $(M,g)$ with torsion $\pm h$.
Denote by $\nabla^g$ the Levi-Civita connection of $g$. We may write $D$ explicitly with respect to the splitting $\T M= TM\oplus T^\star M$, for all $X\in TM$, as follows :
$$
D_X=\left(\begin{array}{cc}
\nabla^g_X&\frac{1}{2}\wedge^2g^{-1}(i_Xh)\\
\frac{1}{2}i_Xh&(\nabla^g_X)^\star
\end{array}\right).
$$

\quad\\
{\bf Definition.} This connection $D$ is called by Gualtieri the generalized Bismut connection associated to $G$.

\quad\\
This connection enable us to give a new carracterisation of twisted generalized kähler manifold.

\quad\\
{\bf Proposition  12 \cite{Gua3}.} If $\J$ be a $G$-orthogonal  generalized almost complex structure, then $(\J,G)$ defines a twisted generalized Kähler structure if and only if :
\begin{enumerate}
\item $D\J=0$ and,

\item the generalized torsion $\TT_D$ is of type $(2,1)+(1,2)$ with respect to $\J$.
\end{enumerate}

\quad\\
{\bf Theorem B.} Let $n\geq 0$. If $(M,G,\I,\J,\K)$ is a twisted generalized hyperkähler $4n$-manifold and  $D$ the generalized Bismut connection, then the generalized almost complex structure $\JJ_D$ on $\ZZ$  is integrable. 

\quad\\
{\bf Proof.} Using proposition 12, we see that both integrability conditions   of theorem A are trivially true. $\square$

\subsection{Levi-Civita connection}

\quad\\
{\bf Theroem  C \cite{Des2}.} Let $n>0$ and let $(M,g,\mathcal Q)$ be a riemannian $4n$-manifold with a generalized almost quaternionic structure such that the Levi-Civita connection $\nabla^g$ preserve $\mathcal Q$, then the generalized almost complex structure $\JJ_{\nabla^g}$ on $\ZZ(\mathcal Q)$  is integrable. 

\quad\\
{\bf Remark.} The case $n=0$ is also treated in \cite{Des2}.

\subsection{Generalized torsion free connection}
 Let $\mathcal Q$ be a generalized almost quaternionic structure on $M$ locally spanned by a generalized almost hypercomplex structure $(\I,\J,\K)$.  Let $G$ be any generalized  metric on $M$  compatible with $\mathcal Q$. In the basis $\T M=C^+\oplus C^-$, an element $u\in\mathcal Q$ is of the form $\left(\begin{array}{cc}u^+&0\\0&u^-\end{array}\right)$. By projection from $C^\pm$ to $TM$, we can consider $u^\pm$ as an almost complex struture on $TM$.  Thus a generalized almost quaternionic  structure gives two almost quaternionic structures namely $Q^\pm=Vect(I^\pm,J^\pm,K^\pm)$. We will note
$$
\begin{array}{cccccc}
f:&Q^-&\lra&\mathcal Q&\lra&Q^+\\
&u^-&\lms&u=\left(\begin{array}{cc}u^+&0\\0&u^-\end{array}\right)&\lms&u^+.
\end{array}
$$
This map induces an algebra isomorphism from $Vect(Id)\oplus Q^-$ to $Vect(Id)\oplus Q^+$.

\quad\\
{\bf Theorem D.} Let $n>1$ and  let $(M,\mathcal Q,G)$ be a generalized almost  quaternionic $4n$-manifold with a generalized metric $G$ compatible with $\mathcal Q$ such that $Q^+=Q^-$. For any generalized torsion free  connection $\nabla$ on $\T M$ compatible with $\mathcal Q$,  $\JJ_\nabla$ is integrable if and only if locally :
\begin{enumerate}
\item[] there exists a generalized hypercomplex structure such that $\nabla \I=\nabla\J=\nabla\K=0$
\item[or] $f=Id$.
\end{enumerate}

 \quad\\
 {\bf Remark.}   From proposition 4 we see that $f=Id$ correspond to $e^{-b}\mathcal Q e^b$ is an almost quaternionic structure; where $b$ is the 2-form associated to $G$.

\quad\\
{\bf Proof of the theorem D.} 
It is clear that if $\nabla\I=\nabla\J=\nabla\K=0$ for a generalized torsion free connection then  both integrability conditions of theorem A are satisfied. As $b$ is not necessarily closed, the integrability of $\mathbb J_\nabla$ when $f=Id$ is not so clear and requires the following lemma.

\quad\\
 {\bf Lemma.} Let $\nabla$ be a generalized torsion free connection  on  $\T M$
compatible with the pseudo-metric. On
the basis $TM\oplus T^\star M$, it takes form: $$
\nabla=\left[\begin{array}{cc} \nabla^1&0\\
L&\nabla^2\end{array}\right]$$ where :

\begin{enumerate}
\item[i)] $\nabla^1$ is a torsion free connection on $TM$ :
$\nabla^1_XY-\nabla^1_YX=[X,Y]$.

\item[ii)]  $\nabla^2$ is the connection on $T^\star M$ induce by
$\nabla^1$ : $$\forall X,Y\in TM,\forall \xi\in T^\star M,\quad
X.<\xi,Y>=<\nabla^2_X\xi,Y>+<\xi,\nabla^1_XY>.$$

\end{enumerate}

\quad\\
{\bf Proof of the lemma.}  $\nabla$ is a generalized torsion free connection
so $\TT(\X,\Y,\ZZ)=0$ $\forall \X,\Y,\ZZ\in\Gamma(\T M)$. In particular for
$(X,\xi,\eta)\in \Gamma(TM)\times \Gamma(T^\star M)\times \Gamma(T^\star M)$  :
$$
\begin{array}{llll}\TT(X,\xi,\eta)=0&\Longleftrightarrow&<\nabla_X\xi-[X,\xi],\eta>=0\\
&\Longleftrightarrow&\nabla_X\xi\textrm{ is a 1-form}
\end{array}
$$
 On the basis $TM\oplus T^\star M$, a torsion free connection $\nabla$ takes form: $$
\nabla=\left[\begin{array}{cc} \nabla^1&0\\
L&\nabla^2\end{array}\right].$$
On the other hand $\forall X,Y,\xi\in \Gamma(TM)\times \Gamma(TM)\times \Gamma(T^\star M)$
we have :
$$
\begin{array}{llll}\TT(X,Y,\xi)=0&\Longleftrightarrow
&<\nabla_XY-\nabla_YX-[X,Y],\xi>=0\\
&\Longleftrightarrow&\nabla^1\textrm{ is torsion free}.
\end{array}
$$
But $\nabla$ is compatible with the pseudometric so :
$$
X.<\xi,Y>=<\nabla^2_X\xi,Y>+<\xi,\nabla^1_XY>.\quad \square
$$
When $f=Id$, any local generalized almost complex structure leaving in $\mathcal Q$ takes form $u=e^b\left[\begin{array}{cc} J&0\\0&-J^\star\end{array}\right]e^{-b}$ for some local almost complex structure $J$. Using the lemma and the fact that  $\nabla$ preserve $\mathcal Q$, a little computation  shows that :
$$
\nabla_X u=e^b\left[\begin{array}{cc} \nabla^1_XJ&0\\0&-(\nabla^1_XJ)^\star)\end{array}\right]e^{-b}.
$$
In particular this means that $\nabla^1$ preserve the almost quaternionic structure $Q=e^{-b}\mathcal Qe^b$. If we note $R^1$ the curvature of the connection $\nabla^1$ and if we differentiate one more time, we have that :
$$\mathcal R(X,Y).u=e^b\left[\begin{array}{cc} R^1(X,Y).J&0\\0&-(R^1(X,Y).J)^\star)\end{array}\right]e^{-b}.
$$
Thus the integrability of $\mathbb J_\nabla$ on $\mathcal Z(\mathcal Q)$ is a consequence of the integrability of $\mathbb J_{\nabla^1}$ on $Z(Q)$ (cf theorem 1 and 2).

\quad\\
It remains to prove the converse. In the basis $C^+\oplus C^-$  the connection $\nabla$  is partitioned into four blocks: $$
\nabla=\left(\begin{array}{ll}
\nabla^+&\nabla^{+-}\\\nabla^{-+}&\nabla^-\end{array}\right).
$$
Since  $\nabla$ is compatible with $\mathcal Q$, for any vector field $X$, we have 
$$
\begin{array}{lclll}
\nabla_X \I&=&&\gamma(X) \J&-\beta(X)\K\\
\nabla_X \J&=&-\gamma(X) \I&&+\alpha(X) \K\\
\nabla_X \K&=&\beta(X) \I&-\alpha(X) \J&
\end{array}
$$
where $\alpha,\beta,\gamma$ are 1-form.   Projecting  on $C^\pm$ we get:
\begin{eqnarray}\left\{\begin{array}{llr}
\nabla^+ I^+&=&\gamma J^+-\beta K^+\\
\nabla^+ J^+&=&-\gamma I^++\alpha K^+\\
\nabla^+ K^+&=&\beta I^+-\alpha J^+\end{array}\right. \textrm{ et
} \left\{\begin{array}{llr}
\nabla^- I^-&=&\gamma J^--\beta K^-\\
\nabla^- J^-&=&-\gamma I^-+\alpha K^-\\
\nabla^- K^-&=&\beta I^--\alpha J^-\end{array}\right.
\end{eqnarray}
Let $G=(g,b)$ be the generalized metric compatible with $\mathcal Q$. Using the lemma and the fact that $C^\pm=\{X+(b\pm g)X\in\T M/X\in TM\}$, we find :
$$
\begin{array}{lcl}
\nabla^+&=&(\pi_+)^{-1}\left(\nabla^1+g^{-1}\Big(\nabla^2(b+g)+L-b(\nabla^1)\Big)\right)\pi\\
\nabla^-&=&(\pi_-)^{-1}\left(\nabla^1-g^{-1}\Big(\nabla^2(b-g)+L-b(\nabla^1)\Big)\right)\pi\\
\nabla^{-+}&=&(\pi_+)^{-1}\left(\nabla^1-g^{-1}\Big(\nabla^2(b+g)+L-b(\nabla^1)\Big)\right)\pi\\
\nabla^{+-}&=&(\pi_-)^{-1}\left(\nabla^1+g^{-1}\Big(\nabla^2(b-g)+L-b(\nabla^1)\Big)\right)\pi\\
\end{array}
$$
But $f$ is an automorphism of $\Sp^2$, so it is a rotation. In a suitable basis we can write:
\begin{eqnarray}
\left\{\begin{array}{clc} I^+=f(I^-)&=&I^-\\
J^+=f(J^-)&=&cJ^-+sK^-\\
K^+=f(K^-)&=&-sJ^-+cK^-
\end{array}
\right.
\textrm{, with $c^2+s^2=1$.}
\end{eqnarray}
On the other hand,  from $\nabla\I\in\mathcal Q$ and $I^+=I^-$, we deduce that $\nabla^{+-}I^+=0$ and
 $\nabla^{-+}I^+=0$, and so $\Big(\nabla^2b+L-b(\nabla^1)\Big)I^+=0$. In particular we decuce that $\nabla^+I^+=\nabla^-I^-$. Now using (2) and (3) we have:
$$
\left\{\begin{array}{cccc} c\gamma+s\beta&=&\gamma\\
s\gamma-c\beta&=&-\beta
\end{array}\right.
\Longleftrightarrow
\left\{\begin{array}{cccc} (c-1)\gamma+s\beta&=&0\\
s\gamma+(1-c)\beta&=&0
\end{array}\right.
$$
On each point either $f=Id$ or $\gamma=\beta=0$. We suppose that $f\neq Id$ and denote by $R^+$ the curvature of the connection $\nabla^+$. In this case we have
$$
\left\{\begin{array}{llc} \nabla^+_X I^+&=&0\\
\nabla^+_X J^+&=&\alpha(X)\; K^+\\
\nabla^+_X      K^+&=&-\alpha(X)\; J^+
\end{array}
\right. \Longrightarrow
\left\{\begin{array}{llc} R^+(X,Y). I^+&=&0\\
R^+(X,Y). J^+&=&\;\;\;d\alpha(X,Y)\; K^+\\
R^+(X,Y). K^+&=&-d\alpha(X,Y)\; J^+
\end{array}
\right.
$$
But if $\JJ_\nabla$ is integrable, then  condition $(C2)$ is true. In particular for $\X=\pi^{-1}_+(X)$ and $\Y=\pi^{-1}_-(Y)$, the projection on $C^+$ gives : 
$$
\left\{\begin{array}{ccc}
R^+(J^+X\w Y+X\w J^-Y).J^+&=&0\\
R^+(K^+X\w Y+X\w K^-Y).K^+&=&0
\end{array}
\right. \Longleftrightarrow \left\{\begin{array}{cccc}
d\alpha(J^+X, Y)&=&-d\alpha(X, J^-Y)&\\
d\alpha(K^+X, Y)&=&-d\alpha(X,K^-Y)
\end{array} \right.
$$
so
$$\begin{array}{ll}
&\left\{\begin{array}{cccc}
d\alpha\Big((cJ^-+sK^-)X, Y\Big)&=&-d\alpha(X, J^-Y)&(L1)\\
d\alpha\Big((-sJ^-+cK^-)X, Y\Big)&=&-d\alpha(X,K^-Y)&(L2)
\end{array}
\right.
\\
\Longrightarrow&d\alpha(J^-X,Y)=-d\alpha\Big(X,(cJ^--sK^-)Y\Big)\quad\quad(L3)=(cL1-sL2)
\end{array}
$$
By symmetry we also have 
$$
\begin{array}{ccccc}
&d\alpha(J^-X, Y)&=&-d\alpha(X, J^+Y)\\
\Longrightarrow&d\alpha(J^-X,Y)&=&-d\alpha\Big(X,(cJ^-+sK^-)Y\Big)&(L4)
\end{array}
$$
Taking $(L3)-(L4)$ we have :
$$0=-2s\;d\alpha(X,K^-Y)
\Longrightarrow  0=d\alpha(X,Y)\quad\forall X,Y\in TM.$$
Thus, $\alpha$ is closed so locally exact : $\alpha=d\theta$ for some locally defined  function $\theta$ on
 $M$, and it is easy to  check that $\Big(\I,\cos \theta\J-\sin\theta\K,\sin\theta\J+\cos\theta\K\Big)$ is a generalized hypercomplex structure such that $
 \nabla\I=\nabla\J=\nabla\K=0.\quad \square$



\quad\\
Guillaume DESCHAMPS,\\
\quad\\
 Université de Brest, UMR 6205,\\ Laboratoire de
Mathématique de
Bretagne Atlantique\\
6 avenue Victor le Gorgeu\\
CS 93837,\\
29238 Brest cedex 3\\
France.

\end{document}